\documentstyle[11pt,leqno,amscd,amssymb,pstricks,latexsym,amsbsy,xypic,mathrsfs,verbatim]{amsart}
\input xy
\xyoption{all}

\newcommand{\scr}
[1]{\mathscr #1}

\newtheorem{Thm}{Theorem}[section]
\newtheorem{Lem}[Thm]{Lemma}
\newtheorem{Cor}[Thm]{Corollary}
\newtheorem{Prop}[Thm]{Proposition}

 \theoremstyle{definition}

\newtheorem{Rem}[Thm]{Remark}

\numberwithin{equation}{section}

\def\lr#1{\langle #1\rangle}
\def\modcat{{\rm mod}}

\def\Hom{{\rm Hom}}
\def\End{{\rm End}}
\def\Ext{{\rm Ext}}

\def\dim{{\rm dim}}
\def\Ker{{\rm Ker}}

\newcommand{\co}{{\mathcal O}}
\newcommand{\vx}{\vec{x}}

\newcommand{\vc}{\vec{c}}

\def\coh{\mbox{\rm coh-}}

\def\az{\alpha}
\def\bz{\beta}
\def\gz{\gamma}

\def\Dim{\mbox{\rm \textbf{dim}}\,}
\def\lz{\lambda}

\def\ez{\epsilon}

\def\Thz{\Theta}

\def\Lz{\Lambda}

\def\bfK{{\bf K}}

\def\bbX{{\mathbb X}}

\def\bbZ{{\mathbb Z}}

\def\bbL{{\mathbb L}}

\def\bbC{{\mathbb C}}
\def\bbP{{\mathbb P}}

\def\A{{\mathcal A}}\def\B{{\mathcal B}}\def\P{{\mathscr{P}}}

\def\cR{{\cal R}}

\def\cE{{\cal E}}
\def\cS{{\cal S}}

\def\cH{{\cal H}}

\def\sC{{\scr C}}

\def\fg{{\frak g}}
\def\cM{{\cal M}}

\def\ra{\rightarrow}

\begin{document}

\title[Exceptional sequences and Drinfeld double Hall algebras]{Exceptional sequences and\\
Drinfeld double Hall algebras}
\author{Shiquan Ruan and Haicheng Zhang }
\address{Yau Mathematical Sciences Center, Tsinghua University,
Beijing 100084,  China.} \email{sqruan@@math.tsinghua.edu.cn,\;\;zhanghai14@@mails.tsinghua.edu.cn }

\thanks{The first author is supported by the Natural Science Foundation of China (Grant No.11471269).}


\begin{abstract}
Let $\A$ be a finitary hereditary abelian category and $D(\A)$ be its reduced Drinfeld double Hall algebra.
By giving explicit formulas in $D(\A)$ for left and right mutations, we show that the subalgebras of $D(\A)$ generated by exceptional sequences are invariant under mutation equivalences. As an application, we obtain that if $\A$ is the category of finite dimensional modules over a finite dimensional hereditary algebra, or the category of coherent sheaves on a weighted projective line, the double composition algebra of $\A$ is generated by any complete exceptional sequence.
Moreover, for the Lie algebra case, we also have paralleled results.
\end{abstract}


\maketitle

\bigskip

\section{Introduction}

Exceptional objects and exceptional sequences have been studied in various contexts. The terminology
{``exceptional"} was first used by Rudakov and his school \cite{Ru} when dealing
with vector bundles on the complex projective plane $\bbP^2$, or more generally for del Pezzo surfaces.
They proposed a useful way to construct exceptional objects from given ones by means of left and right mutations. These mutation operators induce an action of the braid group $B_n$ on the set of exceptional sequences of length $n$.
Bondal and Polishchuk conjectured
in \cite{BP} that the semi-product of the braid group $B_n$ with $\bbZ^n$ acts transitively on the set of complete exceptional sequences for any triangulated category which is generated by an exceptional
sequence of length $n$.

For a finitary category $\A$ over an algebraically closed field, which is a hereditary length category,
Schofield \cite{Schofield} exhibited an algorithm to construct all exceptional objects in $\A$. This algorithm is in fact effective for any field by Ringel \cite{Ringel96}.
By using Schofield's algorithm, Crawley-Boevey \cite{CB} showed that the
braid group acts transitively on the set of complete exceptional sequences in the module category
of a finite dimensional hereditary algebra over an algebraically closed field.
A simplification of the proof and a generalization to hereditary Artin algebras were
given by Ringel \cite{Rex}. Later on,
Meltzer \cite{Mel} proved the transitivity for the category of coherent sheaves over a weighted projective line. More generally, the transitivity also holds for any exceptional curve (cf. \cite{KM}).

Given a finitary hereditary abelian category $\A$, we have the so-called Ringel--Hall algebra $\cH(\A)$ {\cite{Ringel0,Ringel1}}.
For a finite dimensional hereditary algebra $A$ over a finite field, which is of finite representation type, Ringel proved that $\cH(A):=\cH(\modcat A)$ is isomorphic to the positive part 
of the quantized enveloping algebra associated to $A$. For any hereditary algebra $A$, Green \cite{Green} introduced a bialgebra structure on $\cH(A)$
and showed that the composition subalgebra $C(A)$ of $\cH(A)$ generated by simple $A$-modules
gives a realization of the positive part of the corresponding quantized enveloping algebra.
Based on works of Ringel and Green,
Xiao \cite{X97} defined the antipode for the extended Ringel--Hall algebra $\tilde{\cH}(A)$ and proved that $\tilde{\cH}(A)$ has a Hopf algebra structure.
Moreover, he showed that the reduced Drinfeld double Hall algebra of $A$ provideds a realization of the whole quantized enveloping algebra associated to $A$.

The exceptional objects play central roles in the connections of $\cH(A)$ with Lie theory. Indeed, as mentioned above, the composition subalgebra $C(A)$ is defined via simple $A$-modules, which form a complete exceptional sequence when suitably ordered. Moreover, each exceptional module belongs to $C(A)$ for any hereditary algebra $A$ by \cite{ZP,CX}. But in general, an exceptional sequence may not provide building blocks for $C(A)$. The reason is that the simple generators can not be built from a general exceptional sequence by using the Ringel--Hall multiplication.
This problem has been resolved in \cite{CX} by introducing new operators in $\cH(A)$, the so-called left and right derivations, to serve as a kind of {``division"}.
However, the situation is different at the level of categories. Namely,
the subcategories generated by two mutation-equivalent exceptional sequences (in the same orbit under the action of the braid group) always coincide, see for instance \cite{CB}.

The purpose of this paper is to establish a framework rather than $\cH(\A)$ in which mutation-equivalent exceptional sequences play the same role. We intend to deal with the reduced Drinfeld double Hall algebra $D(\A)$, which can be roughly thought as two copies of $\cH(\A)$. In fact, we will show that the subalgebras of $D(\A)$ generated by mutation-equivalent exceptional sequences coincide. This follows from explicit expressions of mutation formulas in $D(\A)$, see Propositions \ref{formula for left mutation} and \ref{formula for right mutation}. These formulas give an accurate recursive algorithm in $D(\A)$ to express each exceptional object in terms of any given complete exceptional sequence. In particular, for a finite dimensional hereditary algebra $A$, we obtain the explicit recursion formulas in $D(A)$ to express each exceptional module via simple modules. We remark that the results in \cite{Sheng} imply that the derived Hall algebra of $\A$, which can be thought as infinite copies of $\cH(\A)$, also fits our desired framework.

The paper is organized as follows. Section 2 gives a brief review on Ringel--Hall algebras and their Drinfeld doubles, as well as the braid group action on the set of exceptional sequences. The main results of this paper are stated in Section 3, but two alternative proofs of Proposition \ref{formula for left mutation} are given in Sections 5 and 6 respectively.
In Section 4 we treat the projective line case as a typical example, and in Section 7, we give the paralleled results in Lie algebra case. The last section includes some applications.

Throughout the paper, let $k$ be a finite field with $q$ elements, and $v=\sqrt{q}$.
Let $\A$ be a finitary (essentially small, Hom-finite) hereditary abelian $k$-category, and we always assume that the endomorphism ring $\End_{\A} X$ for each indecomposable object $X\in\A$ is local. {For each object $X\in\A$ and positive integer $n$, we denote by $nX$ the direct sum of $n$ copies of $X$.} Let $\P$ be the set of isoclasses (isomorphism classes) $[M]$ of objects $M$ in $\A$, and ${\rm ind}\A$ be the complete set of indecomposable objects in $\A$. We choose a representative $V_\alpha\in\alpha$ for each $\alpha\in\P$. Let $\lz_1,\lz_2\in\P$, we denote by $V_{\lz_1+\lz_2}$ the direct sum $V_{\lz_1}\oplus V_{\lz_2}$. Let $K(\A)$ denote the Grothendieck group of $\A$, and we denote by $\mathfrak{r}=\mathfrak{r}_{\A}$ the rank of $K(\A)$. We write $\hat{M}\in K(\A)$ for the class of an object $M\in\A$. For a finite set $X$, we denote by $|X|$ its cardinality. For each $\lz\in \P$ we denote by $a_\lambda$ the cardinality of the automorphism group of $V_\lambda$. Let $A$ be always a finite dimensional hereditary $k$-algebra, and denote by $\modcat A$ the category of finite dimensional (left) $A$-modules. $\bbX$ is always a weighted projective line over $k$, and we denote by $\coh\bbX$ the category of coherent sheaves on $\bbX$.

\section{Preliminaries}

In this section, we review the definitions of Hall algebras, exceptional sequences and quantum binomial coefficients, and collect some necessary known results. Let $\A$ be a finitary hereditary abelian $k$-category.

\subsection{Drinfeld double Hall algebras}

Given objects $A, B\in\A$, let $$\lr{A,B}=\dim_k\Hom_{\A}(A,B)-\dim_k\Ext_{\A}^1(A,B).$$
This descends to a bilinear form $$\lr{-,-}:K(\A)\times K(\A)\longrightarrow\mathbb{Z},$$ known as the \emph{Euler form}. The \emph{symmetric Euler form} $(-,-)$ is given by $(\az,\bz)=\lr{\az,\bz}+\lr{\bz,\az}$ for any $\az,\bz\in K(\A)$.

The \emph{Ringel--Hall algebra} $\mathcal {H}(\A)$ is by definition the $\mathbb{C}$-vector space with the basis $\{u_{\az}|\az\in\P\}$ and with the multiplication given by $u_{\az}u_{\bz}=v^{\lr{\az,\bz}}\sum_{\lambda\in\P}g_{\az\bz}^{\lambda}u_{\lambda}$ for all $\az,\bz\in\P$, where $g_{\az\bz}^{\lambda}$ is the number of subobjects $X$ of $V_{\lambda}$ such that $V_{\lambda}/X$ and $X$ lie in the isoclasses $\az$ and $\bz$, respectively. By Riedtmann-Peng formula \cite{Riedtmann,Peng},
$$g_{\az\bz}^{\lambda}=\frac{|\Ext^1_{\A}(V_{\az},V_{\bz})_{V_\lambda}|}{|\Hom_{\A}(V_{\az},V_{\bz})|}\cdot \frac{a_{\lambda}}{a_{\az}a_{\bz}},$$
where $\Ext^1_{\A}(V_{\az},V_{\bz})_{V_\lambda}$ denotes the subset of $\Ext^1_{\A}(V_{\az},V_{\bz})$ consisting of equivalence classes
of exact sequences in $\A$ of the form $0\ra V_{\bz}\ra V_{\lambda}\ra {V_{\az}}\ra 0.$

Let $\bfK=\bbC[K(\A)]$ be the group algebra of the Grothendieck group $K(\A)$. We denote by $K_{\az}$ (resp. $K_{M}$) the element of $\mathbf{K}$ corresponding to the class $\az\in K(\A)$ (resp. the object $M\in\A$). The \emph{extended Ringel--Hall algebra} is the vector space $\tilde{\mathcal {H}}(\A):=\mathcal {H}(\A)\otimes_{\mathbb{C}}\mathbf{K}$ with the structure of an algebra (containing $\mathcal {H}(\A)$ and $\mathbf{K}$ as subalgebras) by imposing the relations
$$K_{\az}u_M=v^{(\az,\hat{M})}u_MK_{\az}, \text{\ for\ any\ } \az\in K(\A), M\in\A, $$ where and elsewhere we write $u_M$ for $u_{[M]}$ if $M$ is an object in $\A$. Green \cite{Green} introduced a topological bialgebra structure on $\tilde{\mathcal {H}}(\A)$ by defining the comultiplication as follows:
$$\Delta(u_\lambda K_{\gz})
=\sum\limits_{\az,\bz\in\P}v^{\lr{\az,\bz}}\frac{a_{\az}a_{\bz}}{a_{\lambda}}g_{\az\bz}^{\lambda}u_{\az}K_{\bz+\gz}\otimes u_{\bz}K_{\gz}, \text{\ for\ any\ }  \lambda\in\P, \gz\in K(\A).$$
%
%
%
%
%
%
%
%
%
Moreover, there exists the so-called Green's pairing
$(-,-): \tilde{\mathcal {H}}(\A)\times\tilde{\mathcal {H}}(\A)\to \mathbb{C}$
given by the formula $$(u_\lambda K_{\az},u_{\lambda'}K_{\bz})=\frac{\delta_{\lambda,\lambda'}}{a_{\lambda}}v^{-(\az,\bz)}.$$
This pairing is non-degenerate and satisfies that $(ab,c)=(a\otimes b, \Delta(c))$ for any $a,b,c\in\tilde{\mathcal {H}}(\A)$. That is, it is a Hopf pairing.


Now we introduce the reduced Drinfeld double of the topological bialgebra $\tilde{\mathcal {H}}(\A)$.
For this let $\tilde{\mathcal {H}}^{\pm}(\A)$ be two copies of $\tilde{\mathcal {H}}(\A)$ viewed as algebras, and write their basis 
as $\{u_\lambda^{\pm}K_{\gz}|\lambda\in\P, \gz\in K(\A)\}$, and the comultiplications are given by
$$\Delta(u_\lambda^{\pm} K_{\gz})
=\sum\limits_{\az,\bz\in\P}v^{\lr{\az,\bz}}\frac{a_{\az}a_{\bz}}{a_{\lambda}}g_{\az\bz}^{\lambda}u_{\az}^{\pm}K_{^{\pm}\bz+\gz}\otimes u_{\bz}^{\pm}K_{\gz}. $$

The \emph{Drinfeld double} of the topological bialgebra $\tilde{\mathcal {H}}(\A)$ with respect
to the Green's pairing $(-,-)$ is the associative algebra $\tilde{D}(\A)$, defined as the free product of
algebras $\tilde{\mathcal {H}}^+(\A)$ and $\tilde{\mathcal {H}}^-(\A)$ subject to the following relations $D(a,b)$ for all $a,b\in\tilde{\mathcal {H}}(\A)$:
\begin{equation}\label{Drinfeld relation}\sum(a_1,b_2)b_1^+\ast a_2^-=\sum(a_2,b_1)a_1^-\ast b_2^+,\end{equation}
where we write $\Delta(c^{\pm})=\sum c_1^{\pm}\otimes c_2^{\pm}$ in the Sweedler's notation for $c\in\tilde{\mathcal {H}}(\A)$ (see for example \cite{BurSch}). 
The \emph{reduced Drinfeld double} $D(\A)$ is the quotient of $\tilde{D}(\A)$ by the two-sided ideal
$$I=\langle {K_{\az}\otimes1-1\otimes K_{-\az} } , \az\in K(\A)\rangle.$$

For each object $M\in\A$ and $t\in\mathbb{Z}$, we sometimes also write $K_M^{\pm t}$ for $K_{\pm t\hat{M}}$.
We have the following fundamental result of the reduced Drinfeld double Hall algebras for finitary hereditary abelian categories due to Cramer.


\begin{Thm}{\rm(\cite{cram})}\label{cramer}
Let $\A$ and $\B$ be two finitary hereditary abelian $k$-categories such that there
exists an equivalence of triangulated categories $F: D^b(\A)\cong D^b(\B)$.
Then there is an algebra isomorphism $\Phi: D(\A)\cong D(\B)$ uniquely determined by the following property.
For any $\az\in K(\A)$, we have $\Phi(K_{\az})=K_{F(\az)}$. Next, for any object $M\in\mathcal{A}$ such that $F(M)=N[-n]$ with $N\in\mathcal{B}$ and $n\in\bbZ$, we have
$$\Phi(u_{M}^{\pm})=v^{-n\langle N,N\rangle}u_{N}^{\pm \bar{n}}K_{N}^{\pm n},$$
where $\bar{n}=+$ (resp. $-$) if $n$ is even (resp. odd).
\end{Thm}

\subsection{Exceptional sequences}

An object $V_{\az}\in\A$ is called \emph{exceptional} if $\Ext_{\A}^1(V_{\az},V_{\az})$ vanishes and $\End_{\A}V_{\az}$ is a skew field. In our case, it is a field. A pair $(V_{\az},V_{\bz})$ of exceptional objects in $\A$ is called \emph{exceptional} provided $$\Hom_{\A}(V_{\bz},V_{\az})=\Ext_{\A}^1(V_{\bz},V_{\az})=0.$$ A sequence $(V_{\az_1},V_{\az_2},\cdots,V_{\az_r})$ of exceptional objects in $\A$ is called \emph{exceptional} provided any pair $(V_{\az_i},V_{\az_j})$ with $1\leq i<j\leq r$ is exceptional.
It is said to be \emph{complete} if $r=\mathfrak{r}_{\A}$, the rank of the Grothendieck group $K(\A)$.
For any exceptional pair $(V_{\az},V_{\bz})$, let $\mathscr{C}$$(V_{\az},V_{\bz})$ be
the smallest full subcategory of $\A$ which contains the objects $V_{\az}, V_{\bz}$, and
is closed under kernels, cokernels and extensions. As an exact extension-closed full subcategory of $\A$, $\sC(V_{\az},V_{\bz})$ is itself abelian and hereditary, which is 
derived equivalent to the module category $\modcat \Lz$, where and elsewhere $\Lz$ is a finite dimensional
hereditary $k$-algebra with two isoclasses of
simple modules {(\cite{HR1,HR2,LM}).}
%
%
%
%
%
%
%
%

According to {\cite{Bondal,CB,CX,Rex}}, for any exceptional pair $(V_{\az},V_{\bz})$, there exist unique objects $L(\az,\bz)$ and $R(\az,\bz)$ with the property that $(L(\az,\bz),V_{\az})$ and $(V_{\bz},R(\az,\bz))$ are both exceptional pairs in $\A$. The objects $L(\az,\bz)$ and $R(\az,\bz)$ are called the \emph{left mutation} and \emph{right mutation} of $(V_{\az},V_{\bz})$, respectively.
Moreover, if $\cE=(V_{\az_1},V_{\az_2},\cdots,V_{\az_r})$ is an exceptional sequence, for $1\leq i<r$, we define
$$L_i\cE:=(V_{\az_1},\cdots,V_{\az_{i-1}},L(\az_i,\az_{i+1}),V_{\az_i},V_{\az_{i+2}},\cdots,V_{\az_r}),$$
$$R_i\cE:=(V_{\az_1},\cdots,V_{\az_{i-1}},V_{\az_{i+1}},R(\az_i,\az_{i+1}),V_{\az_{i+2}},\cdots,V_{\az_r}).$$
Then $L_i\cE$ and $R_i\cE$ are both exceptional sequences.

Recall that the braid group $\mathcal {B}$$_r$ in $r-1$ generators $\sigma_1,\sigma_2,\cdots,\sigma_{r-1}$ is the group with these generators and the relations $\sigma_{i}\sigma_{i+1}\sigma_{i}=\sigma_{i+1}\sigma_{i}\sigma_{i+1}$ for all $1\leq i<r-1$ and $\sigma_{i}\sigma_{j}=\sigma_{j}\sigma_{i}$ for $j\geq i+2$. The braid group $\mathcal {B}$$_r$ acts on the set of all exceptional sequences of length $r$ via $\sigma_i\cE:=L_i\cE$ and $\sigma_i^{-1}\cE:=R_i\cE$.
Two exceptional sequences of length $r$ are said to be \emph{mutation equivalent} if they belong to the same $\mathcal {B}$$_r$-orbit.

\subsection{Quantum binomial coefficients}

In Hall algebras, the following notations and relations are frequently-used. Let $i,l$ be two nonnegative integers and $0\leq i\leq l$, set
\begin{equation*}
[l]=\frac{v^l-v^{-l}}{v-v^{-1}},~[l]!=\prod_{t=1}^l[t],~\left[\begin{smallmatrix}
l\\i
\end{smallmatrix}\right]=\frac{[l]!}{[i]![l-i]!},
\end{equation*}
and

\begin{equation*}
|l]=\frac{{q}^l-1}{{q}-1},~|l]!=\prod_{t=1}^l|t],~\left|\begin{smallmatrix}
l\\i
\end{smallmatrix}\right]=\frac{|l]!}{|i]!|l-i]!}.
\end{equation*}
Then $|l]=v^{l-1}[l],~|l]!=v^{\left(\begin{smallmatrix}
l\\2
\end{smallmatrix}\right)}[l]!,~\left|\begin{smallmatrix}
l\\i
\end{smallmatrix}\right]=v^{i(l-i)}\left[\begin{smallmatrix}
l\\i
\end{smallmatrix}\right].$ For a polynomial $f\in$$\mathbb{Z}$$[v,v^{-1}]$ and a positive integer $d$, we denote by $f(v)_d$ the polynomial obtained from $f$ by replacing $v$ by $v^d$.

The following lemma is well-known (see for example \cite{R96}).
\begin{Lem}\label{zuhe}
For any integer $l>0$,
\begin{equation}
\sum_{i=0}^l(-1)^iv^{i(i-1)}\left|\begin{smallmatrix}
l\\i
\end{smallmatrix}\right]=0.
\end{equation}
\end{Lem}

For any exceptional object $V_{\lz}\in\A$, set $u_{\lz}^{\pm(t)}=(1/{[t]^!_{\epsilon(\lz)}})u_{\lz}^{\pm t}$ in the reduced Drinfeld double Hall algebra $D(\A)$, where $\epsilon(\lz)=\dim_k\End_{\A}V_{\lz}$.
We have the identity $u_{\lz}^{\pm(t)}=(v^{\epsilon(\lz)})^{t(t-1)}u_{t\lz}^{\pm}$.

\section{Main results}

In this section we state our main results of this paper. For any exceptional sequence $\cE=(V_{\az_1},V_{\az_2},\cdots,V_{\az_{r}})$ in $\A$, we denote by $D_{\cE}(\A)$ the subalgebra of $D(\A)$ generated by the elements $\{u_{\az_i}^\pm, K_{\az_i}|1\leq i\leq r\}$, and call $D_{\cE}(\A)$ the subalgebra generated by $\cE$. 

\begin{Thm}\label{Main theorem}
Let $\A$ be a finitary hereditary abelian $k$-category. If two exceptional sequences $\cE_1$ and $\cE_2$ are mutation equivalent, then $D_{\cE_1}(\A)=D_{\cE_2}(\A)$.
\end{Thm}

\begin{pf}
Let $r$ be the length of $\epsilon_1$ and $\epsilon_2$. By induction, we assume that $\cE_2$ is obtained from $\cE_1$ by one step left mutation (resp. one step right mutation), i.e., $\cE_2=L_i(\cE_1)$ (resp. $\cE_2=R_i(\cE_1)$) for some $1\leq i\leq r-1$. Moreover, by the definitions of $L_i$ and $R_i$, we know that $\cE_1$ and $\cE_2$ only differ in two components. Thus we assume that $\cE_1$ and $\cE_2$ have length $r=2$ without loss of generality. In this case, by \cite{HR1,HR2,LM} we know that the category $\mathscr{C}$$(\cE_i)$ is derived equivalent to $\modcat \Lz$ for some hereditary $k$-algebra $\Lz$ with two isoclasses of simple modules for $i=1,2$. By Cramer's Theorem \ref{cramer}, which states that the Drinfeld double Hall algebra is derived invariant, we further assume that $\mathscr{C}$$(\cE_1)$ is just $\modcat \Lz$. Now the result follows from Proposition \ref{formula for left mutation} (\emph{resp}. Proposition \ref{formula for right mutation}), which gives an explicit expression of $u_{\gz}^{\pm}$ for the left mutation $L(\az,\bz)=V_{\gz}$ (\emph{resp}. of $u_{\lz}^{\pm}$ for the right mutation $R(\az,\bz)=V_{\lz}$) in $u_{\az}^{\pm}, u_{\bz}^{\pm}, K_{\pm\az}, K_{\pm\bz}$ for any exceptional pair $(V_{\az},V_{\bz})$ in $\modcat \Lambda$.
\end{pf}

Let $(V_{\az},V_{\bz})$ be an exceptional pair in $\modcat \Lz$. Denote by $\epsilon(\az)=\lr{\az,\az}, \epsilon(\bz)=\lr{\bz,\bz},$ and $$n(\az,\bz)=\frac{\lr{\az,\bz}}{\lr{\az,\az}},\ m(\az,\bz)=\frac{\lr{\az,\bz}}{\lr{\bz,\bz}},\ n=|n(\az,\bz)|, \ m=|m(\az,\bz)|.$$
The following two propositions provide explicit formulas in the reduced Drinfeld double Hall algebra for the left mutation and right mutation respectively.

\begin{Prop}\label{formula for left mutation}
Let $(V_{\az},V_{\bz})$ be an exceptional pair in $\modcat \Lambda$ with the left mutation $L(\az,\bz)=V_{\gz}$.
\begin{itemize}
\item[(i)] If $\lr{\az,\bz}\leq 0$, then
$$u_{\gz}^{\pm}=\sum_{j=0}^n(-1)^j(v^{\epsilon(\az)})^{n-j}u_{\az}^{\pm(n-j)}u_{\bz}^{\pm}
u_{\az}^{\pm(j)};$$
\item[(ii)] If $n\Dim V_{\az}>\Dim V_{\bz}$, then $$u_{\gz}^{\pm}=v^{\epsilon(\bz)}K_{\pm\bz}\sum_{j=0}^n(-1)^j(v^{\epsilon(\az)})^{-(n-1)(n-j)}u_{\az}^{\pm(n-j)}u_{\bz}^{\mp}
u_{\az}^{\pm(j)};$$
\item[(iii)] If $0<n\Dim V_{\az}<\Dim V_{\bz}$, then
$$u_{\gz}^{\pm}=K_{\mp n\az}\sum_{j=0}^n(-1)^j(v^{\epsilon(\az)})^{-(n-1)(n-j)}u_{\az}^{\mp(n-j)}u_{\bz}^{\pm}
u_{\az}^{\mp(j)}.$$
\end{itemize}
\end{Prop}

\begin{Prop}\label{formula for right mutation}
Let $(V_{\az},V_{\bz})$ be an exceptional pair in $\modcat \Lambda$ with the right mutation $R(\az,\bz)=V_{\lz}$.
\begin{itemize}
\item[(i)] If $\lr{\az,\bz}\leq 0$, then
$$u_{\lz}^{\pm}=\sum_{j=0}^m(-1)^j(v^{\epsilon(\bz)})^{m-j}u_{\bz}^{\pm(j)}u_{\az}^{\pm}
u_{\bz}^{\pm(m-j)};$$
\item[(ii)] If $m\Dim V_{\bz}>\Dim V_{\az}$, then $$u_{\lz}^{\pm}=v^{\epsilon(\az)}\sum_{j=0}^m(-1)^j(v^{\epsilon(\bz)})^{-(m-1)(m-j)}u_{\bz}^{\pm(j)}u_{\az}^{\mp}
u_{\bz}^{\pm(m-j)}K_{\mp\az};$$
\item[(iii)] If $0<m\Dim V_{\bz}<\Dim V_{\az}$, then
$$u_{\lz}^{\pm}=\sum_{j=0}^m(-1)^j(v^{\epsilon(\bz)})^{-(m-1)(m-j)}u_{\bz}^{\mp(j)}u_{\az}^{\pm}
u_{\bz}^{\mp(m-j)}K_{\pm m\bz}.$$
\end{itemize}
\end{Prop}

Let $D_{\mathscr{E}}(\A)$ be the subalgebra of $D(\A)$ generated by all $K_{\az},\az\in K(\A)$, and all $u_{\lz}^{\pm}$ corresponding to exceptional objects $V_{\lz}$ in $\A$. {We call $D_{\mathscr{E}}(\A)$ the \emph{double composition algebra} of $\A$.}

\begin{Prop}\label{main cor}
Let $\mathcal{A}$ and $\mathcal{B}$ be two finitary hereditary abelian $k$-categories. If $\mathcal{A}$ and $\mathcal{B}$ are derived
equivalent, then $D_{\mathscr{E}}(\A)\cong D_{\mathscr{E}}(\mathcal{B})$.
\end{Prop}

\begin{pf}
By Cramer's Theorem \ref{cramer}, the derived equivalence between $\mathcal{A}$ and $\mathcal{B}$ implies an isomorphism between their Drinfeld double Hall algebras. Moreover, note that the derived equivalence preserves exceptional objects in both categories. Thus $D_{\mathscr{E}}(\A)\cong D_{\mathscr{E}}(\mathcal{B})$. This finishes the proof.
\end{pf}

\begin{Thm}\label{add conditions for A}
Let $\A$ be a finitary hereditary abelian $k$-category. Assume that $\A$ satisfies the following conditions:
\begin{itemize}
\item[(i)] there exists a complete exceptional sequence in $\A$;
\item[(ii)] the action of the braid group $\mathcal {B}$$_{\mathfrak{r}}$ on the set of complete exceptional sequences in $\A$ is transitive;
\item[(iii)] for any exceptional object $E\in\A$, there exists a complete exceptional sequence $\cE'$ in $\A$ such that $E\in\cE'$.
\end{itemize}
Then for any complete exceptional sequence $\cE$ in $\A$, $D_{\mathscr{E}}(\A)=D_{\cE}(\A)$.
\end{Thm}

\begin{pf}
Obviously, by definition we have an inclusion $D_{\cE}(\A)\subseteq D_{\mathscr{E}}(\A)$. To finish the proof we only need to show that for each exceptional object $E$ in $\A$, the generator $u_{E}^{\pm}$ of $D_{\mathscr{E}}(\A)$ belongs to $D_{\cE}(\A)$. For this, choose a complete exceptional sequence $\cE'$ in $\A$ such that $E\in\cE'$. By transitivity we obtain that $D_{\cE}(\A)=D_{\cE'}(\A)$ according to Theorem \ref{Main theorem}. This finishes the proof.
\end{pf}

\begin{Cor}\label{main examples}
Let $\A=\modcat A$ for some finite dimensional hereditary algebra $A$ or $\A=\coh\bbX$ for some weighted projective line $\bbX$. Then for any complete exceptional sequence $\cE$ in $\A$, $D_{\mathscr{E}}(\A)=D_{\cE}(\A)$.
\end{Cor}

\begin{pf}
By \cite{CB,Rex} and \cite{KM} we know that $\A$ satisfies the conditions in Theorem \ref{add conditions for A} in both cases. This finishes the proof.
\end{pf}


{\begin{Rem}(1) In the case $\A=\modcat A$, the simple $A$-modules form a complete exceptional sequence when suitably ordered. Hence the double composition algebra of $\A$ can be defined via simple modules, which in fact has been adopted as its original definition, see for example \cite{Ringel1, X97}.
Moreover, in this case, the above result is closely related to \cite[Corollary 5.3]{CX}. For this we only need to observe that, under the natural embedding $$\cH(\A)\hookrightarrow D(\A); u_{\az}\mapsto u^+_{\az}, \az\in\P,$$ the left and right derivations $_{i}\delta, \delta_{i}$ of $\cH(\A)$ associated to each simple $A$-module $S_i$ are related to the adjoint operator $[u_{i}^-, -]$ (see for example \cite[Proposition 3.1.6]{Lus}).

(2) In the case $\A=\coh\bbX$, the double composition algebra $D_{\mathscr{E}}(\A)$ appears in different contexts, which will be discussed in detail in Subsection \ref{composition algebra of coh subsection}.
\end{Rem}}

\section{Typical example}

In this section, we consider $\A$ as the category $\coh\bbP^1$ of coherent sheaves on the projective line $\mathbb{P}$$^1$ over $k$. The Ringel--Hall of $\coh\bbP^1$ has been widely studied, which are closely related to the quantum affine algebra $U_{v}(\hat{sl_2})$, see for example \cite{Ka, BK, Sch2006, BurSch}. This is our original example for Theorem \ref{Main theorem}, and the content in this section itself is of interest.

For the category $\coh\bbP^1$, it is well-known that the exceptional objects coincide with the indecomposable vector bundles, which are all line bundles by \cite{Grothendieck} and hence have the form $\co(i), i\in\mathbb{Z}$. Moreover, each exceptional pair in $\coh\bbP^1$ has the form $(\co(i),\co(i+1)), i\in\mathbb{Z}$. In this case $n=n(\co(i),\co(i+1))=2$ and $m=m(\co(i),\co(i+1))=2$.
The left (resp. right) mutation of $(\co(i),\co(i+1))$ is $\co(i-1)$ (resp. $\co(i+2)$), which comes from the following Auslander--Reiten sequence for $j=i$ (resp. $j=i+1$):
$$0\longrightarrow \co(j-1)\longrightarrow \co(j)\oplus\co(j)\longrightarrow \co(j+1)\longrightarrow 0.$$

We will give explicit expressions in the reduced Drinfeld double Hall algebra of $\coh\bbP^1$ for the elements $u_{\co(i-1)}^{\pm}$ and $u_{\co(i+2)}^{\pm}$ corresponding to the left and right mutations of the exceptional pair $(\co(i),\co(i+1))$ respectively.


Recall from \cite{Sch2006} that for any $j\in\mathbb{Z}$,
\begin{equation}\label{yu1}\Delta(u_{\co(j)}^\pm)=u_{\co(j)}^{\pm}\otimes 1+K_{{\co}(j)}^{\pm}\otimes u_{\co(j)}^{\pm}+\sum_{k>0}\Thz_k^{\pm}K_{{\co}(j-k)}^{\pm}\otimes u_{\co(j-k)}^{\pm}.\end{equation}
The elements $\Thz_k^{\pm}, k\geq 0$ play important roles in the study of $D(\coh\bbP^1)$, for its concrete definition we refer to \cite[Example 4.12]{Sch2006}. We emphasize that $\Thz_1^{\pm}=(v-v^{-1})\sum\limits_{S\in\cS_1}u_{S}^{\pm}$, where $\cS_1$ denotes the set of torsion sheaves of degree one, and $|\cS_1|=q+1$. Let {$\delta= \widehat{{\co}(1)}-\widehat{{\co}}$}. Then

\begin{equation}\label{yu2}\Delta(\Thz_1^\pm)=\Thz_1^\pm\otimes1+K_{\pm\delta}\otimes\Thz_1^\pm~~\mbox{and}~~ (\Thz_1^+,\Thz_1^-)=q-q^{-1}.\end{equation}

\begin{Prop}\label{left formula for P1} For any $i\in\mathbb{Z}$,
$$u_{\co(i-1)}^{\pm}=K_{{\co}(i+1)}^{\pm}\sum_{j=0}^2(-1)^jv^{j-1}u_{\co(i)}^{\pm(2-j)}
u_{\co(i+1)}^{\mp}u_{\co(i)}^{\pm(j)}.$$\end{Prop}

\begin{pf} By duality, we only show the formula for $u_{\co(i-1)}^{+}$.
Using $(\ref{yu1})$ and $(\ref{yu2})$, we obtain that
\begin{equation*}\begin{split}
u_{\co(i)}^+u_{\co(i+1)}^-&=u_{\co(i+1)}^-u_{\co(i)}^++\Thz_1^-K_{{\co}(i)}^{-}(u_{\co(i)}^+,u_{\co(i)}^-)\\&=
u_{\co(i+1)}^-u_{\co(i)}^++(q-1)^{-1}\Thz_1^-K_{{\co}(i)}^{-},
\end{split}\end{equation*} and
\begin{equation*}\begin{split}
u_{\co(i)}^+\Thz_1^-&=\Thz_1^-u_{\co(i)}^++K_{-\delta}u_{\co(i-1)}^+(\Thz_1^+K_{{\co}(i-1)},\Thz_1^-)\\&=
\Thz_1^-u_{\co(i)}^++(q-q^{-1})K_{-\delta}u_{\co(i-1)}^+.
\end{split}\end{equation*}
Hence, $\Thz_1^-=(q-1)[u_{\co(i)}^+, u_{\co(i+1)}^-]K_{{\co}(i)}$ and
$u_{\co(i-1)}^+=(q-q^{-1})^{-1}K_{\delta}[u_{\co(i)}^+, \Thz_1^-].$
Observing that $K_{\delta}$ is central and $K_{{\co}(i)}K_{\delta}=K_{{\co}(i+1)}$, we obtain that
$$u_{\co(i-1)}^+=K_{{\co}(i+1)}
(v^{-1}u_{\co(i)}^{+(2)}u_{\co(i+1)}^--u_{\co(i)}^+
u_{\co(i+1)}^-u_{\co(i)}^++vu_{\co(i+1)}^-u_{\co(i)}^{+(2)}).$$

%
\end{pf}

\begin{Prop} For any $i\in\mathbb{Z}$,
$$u_{\co(i+2)}^{\pm}=\sum_{j=0}^2(-1)^jv^{j-1}u_{\co(i+1)}^{\pm(j)}
u_{\co(i)}^{\mp}u_{\co(i+1)}^{\pm(2-j)}K_{{\co}(i)}^{\mp}.$$\end{Prop}

\begin{pf}
By duality, we only show the formula for $u_{\co(i+2)}^{+}$. Using similar arguments as those in the proof of Proposition \ref{left formula for P1}, the element $\Thz_1^+$ can be expressed by the elements in $D(\coh\bbP^1)$ corresponding to the exceptional pair $(\co(i), \co(i+1))$ as follows:
$$\Thz_1^+=(q-1)[u_{\co(i)}^-,u_{\co(i+1)}^+]K_{{\co}(i)}^-.$$
On the other hand, by definition $\Thz_1^+=(v-v^{-1})\sum\limits_{S\in\cS_1}u_{S}^+$, and for any $S\in\cS_1$,
$$g_{\co(i),S}^{\co(i)\oplus S}=g_{S,\co(i+1)}^{\co(i+2)}=1,\ \ g_{S,\co(i+1)}^{S\oplus \co(i+1)}=q.$$
Hence,
\begin{equation*}\begin{split}
\Thz_1^+u_{\co(i+1)}^+&=(v-v^{-1})\sum\limits_{S\in\cS_1}u_{S}^+u_{\co(i+1)}^+\\
&=(v-v^{-1})\sum\limits_{S\in\cS_1}v^{\lr{S,\co(i+1)}}(qu_{S\oplus \co(i+1)}^++u_{\co(i+2)}^+)\\
&=(v-v^{-1})\sum\limits_{S\in\cS_1}(vu_{S\oplus \co(i+1)}^++v^{-1}u_{\co(i+2)}^+)\\
&=(v-v^{-1})\sum\limits_{S\in\cS_1}u_{\co(i+1)}^+u_{S}^++(1-q^{-1})\sum\limits_{S\in\cS_1}
u_{\co(i+2)}^+\\
&=u_{\co(i+1)}^+\Thz_1^++(1-q^{-1})(1+q)u_{\co(i+2)}^+.
\end{split}\end{equation*}
Therefore, \begin{equation*}\begin{split}
u_{\co(i+2)}^+&=(1-q^{-1})^{-1}(1+q)^{-1}[\Thz_1^+, u_{\co(i+1)}^+]\\
&=q(1+q)^{-1}[[u_{\co(i)}^-,u_{\co(i+1)}^+]K_{{\co}(i)}^-,u_{\co(i+1)}^+]\\
&=(vu_{\co(i+1)}^{+(2)}u_{\co(i)}^--u_{\co(i+1)}^+u_{\co(i)}^-u_{\co(i+1)}^++v^{-1}u_{\co(i)}^-u_{\co(i+1)}^{+(2)})K_{{\co}(i)}^-.
\end{split}\end{equation*}
\end{pf}

\section{The first proof for Proposition \ref{formula for left mutation}}

In this section, we present a fundamental proof for Proposition \ref{formula for left mutation} by using explicit comultiplicatin formulas in the reduced Drinfeld double Hall algebra $D(\Lz):=D(\modcat \Lambda)$ of the module category $\modcat \Lambda$.

For the exceptional pair $(V_{\az},V_{\bz})$ in $\modcat \Lambda$, if $\lr{\az,\bz}\leq 0$, then $(V_{\az},V_{\bz})$ are two simple modules in $\modcat \Lambda$, where $V_{\bz}$ is the simple projective module and $V_{\az}$ is the simple injective module. In this case, the formula in Proposition \ref{formula for left mutation} (i) is well-known, see for example \cite[Proposition 4.3.3]{CX}. Hence we only need to prove the statements of Proposition \ref{formula for left mutation} (ii) and (iii) in the following, in which case $n=n(\az,\bz)$.

%
%
%

\subsection{The case of $n\Dim V_{\az}>\Dim V_{\bz}$}

That is, $(V_{\az},V_{\bz})$ are the slice modules in the preprojective or preinjective component of $\modcat \Lambda$.
The left mutation $V_{\gamma}$ of $(V_{\az},V_{\bz})$ is given by the following Auslander--Reiten sequence by \cite{CX}
$$0\to V_{\gamma}\to V_{n\az}\to V_{\bz}\to 0.$$

%

We have the following result:
\begin{Lem}{\rm (\cite{CX})}\label{lem4.2}
If $f: V_{n\az}\to V_{\bz}$ is an epimorphism, then $\Ker f\cong V_{\lambda_i}\oplus V_{i\az}$ with $0\leq i\leq n-1$, where $V_{\lambda_0}\cong V_{\gz}$. Moreover, each $V_{\lambda_i}$ is a submodule of $V_{\gz}$ and has no direct summands which are isomorphic to $V_{\az}$.
\end{Lem}

We will give the proof of Proposition \ref{formula for left mutation} (ii) by considering whether $(V_{\az},V_{\bz})$ lies in the preinjective component or in the preprojective. We only give the proof for the result of $u_{\gz}^+$, since the formula for $u_{\gz}^-$ can be obtained dually.

%

\subsubsection{$(V_{\az},V_{\bz})$ lies in the preinjective component}



In this case, for any $1\leq j\leq n$, each non-zero morphism from $V_{j\az}$ to $V_{\bz}$ is an epimorphism, whose kernel has similar description as above:

\begin{Lem}\label{ker}
Let $1\leq j\leq n$ and $f:V_{j\az}\to V_{\bz}$ be an epimorphism, then $\Ker(f)\cong V_{\lambda_{j,i}}\oplus V_{i\az}$ for some $0\leq i\leq j-1$, where $V_{\lambda_{j,i}}$ contains no direct summands which are isomorphic to $V_{\az}$, and it fits into the following exact sequence $$\xymatrix{0\ar[r]&V_{\lz_{j,i}}\ar[r]&V_{(j-i)\az}\ar[r]&V_{\bz}\ar[r]&0.}$$
Moreover, $[\Ext_{\Lz}^1(V_{\bz},V_{\lz_{j,i}}):\End_{\Lz}V_{\bz}]=1$.
\end{Lem}

\begin{pf}
For the epimorphism $f:V_{j\az}\to V_{\bz}$, we have the exact sequence
$$\xymatrix{0\ar[r]& \Ker f\oplus V_{(n-j)\az}\ar[r]& V_{j\az}\oplus V_{(n-j)\az}\ar[r]^-{f\choose 0}& V_{\bz}\ar[r]& 0.}$$ Then in a similar way to \cite[Lemma 4.1.2]{CX}, we can prove the first statement.
Since $\Hom_{\Lz}(V_{\bz},V_{\lz_{j,i}})=0$, we obtain that
$\lr{V_{\bz},V_{\lz_{j,i}}}=-\dim_k\Ext_{\Lz}^1(V_{\bz},V_{\lz_{j,i}}).$ On the other hand, $\lr{V_{\bz},V_{\lz_{j,i}}}=\lr{\bz,(j-i)\az-\bz}=-\epsilon(\bz)=-\dim_k\End_{\Lz}(V_{\bz})$.
Hence, $[\Ext_{\Lz}^1(V_{\bz},V_{\lz_{j,i}}):\End_{\Lz}V_{\bz}]=1$.
\end{pf}

\begin{Lem}\label{yin1}
Let $0\leq j\leq n-1$. Then
$$u_{(n-j)\az}^+u_{\bz}^-=u_{\bz}^-u_{(n-j)\az}^++\sum_{i=0}^{n-j-1}\sum_{\lz_{n-j,i}:g_{\bz,\lz_{n-j,i}}^{(n-j-i)\az}\neq0}v^{-\epsilon(\bz)}K_{-\bz}
u_{\lambda_{n-j,i}+i\az}^+.$$
\end{Lem}

\begin{pf}
Note that $$\Delta(u_{\bz}^-)=u_{\bz}^-\otimes 1+K_{-\bz}\otimes u_{\bz}^-+\mbox{else}$$ and
$$\Delta(u_{(n-j)\az}^+)=u_{(n-j)\az}^+\otimes 1+K_{(n-j)\az}\otimes u_{(n-j)\az}^++\sum v^{-\epsilon(\bz)}a_{n-j}u_{\bz}^+K_{(n-j)\az-\bz}\otimes u_{k_{n-j}}^++\mbox{else},$$
where the sum on the right-hand side is taken over all the isoclasses $k_{n-j}$'s satisfying $g_{\bz, k_{n-j}}^{(n-j)\az}\neq 0$, and by Lemma \ref{ker}, $\lr{{\bz},k_{n-j}}=\lr{\bz,(n-j)\az-\bz}=-\epsilon(\bz)$ and
\begin{equation*}\begin{split}a_{n-j}&=\frac{|\Ext_{\Lz}^1(V_{\bz},V_{\lz_{n-j,i}}\oplus V_{i\az})_{V_{(n-j)\az}}|}{|\Hom_{\Lz}(V_{\bz},V_{\lz_{n-j,i}}\oplus V_{i\az})|}
\\&=\frac{|\Ext_{\Lz}^1(V_{\bz},V_{\lz_{n-j,i}})_{V_{(n-j-i)\az}}|}{|\Hom_{\Lz}(V_{\bz},V_{\lz_{n-j,i}})|}\\&=q^{\epsilon(\bz)}-1.\end{split}\end{equation*}
Thus, by $(\ref{Drinfeld relation})$ we obtain the desired identity.
\end{pf}

\begin{Lem}\label{yin2}
Let $0\leq j\leq n-1$. Then
\begin{equation*}\begin{split}u_{(n-j)\az}^+u_{\bz}^-u_{j\az}^+&=
(v^{\epsilon(\az)})^{j(n-j)}\left|\begin{smallmatrix}
n\\j
\end{smallmatrix}\right]u_{\bz}^-u_{n\az}^++\\&(v^{\epsilon(\az)})^{j(n-j)}
\sum_{i=0}^{n-j-1}\sum_{\lz_{n-j,i}:g_{\bz,\lz_{n-j,i}}^{(n-j-i)\az}\neq0}v^{-\epsilon(\bz)}K_{-\bz}
\left|\begin{smallmatrix}
i+j\\j
\end{smallmatrix}\right]u_{\lz_{n-j,i}+(i+j)\az}^+.\end{split}\end{equation*}
\end{Lem}

\begin{pf}
Observing that $$u_{(n-j)\az}^+u_{j\az}^+=(v^{\epsilon(\az)})^{j(n-j)}g_{(n-j)\az,j\az}^{n\az}u_{n\az}^+=(v^{\epsilon(\az)})^{j(n-j)}
\left|\begin{smallmatrix}
n\\j
\end{smallmatrix}\right]u_{n\az}^+$$
and \begin{equation*}\begin{split}u_{\lz_{n-j,i}+i\az}^+u_{j\az}^+&=v^{\lr{(n-j)\az-\bz,j\az}}g_{\lz_{n-j,i}+i\az,j\az}^{\lz_{n-j,i}+(i+j)\az}
u_{\lz_{n-j,i}+(i+j)\az}^+\\&
=(v^{\epsilon(\az)})^{j(n-j)}g_{i\az,j\az}^{(i+j)\az}u_{\lz_{n-j,i}+(i+j)\az}^+\\&
=(v^{\epsilon(\az)})^{j(n-j)}\left|\begin{smallmatrix}
i+j\\j
\end{smallmatrix}\right]u_{\lz_{n-j,i}+(i+j)\az}^+,\end{split}\end{equation*}
we finish the proof by Lemma $\ref{yin1}$.
\end{pf}

$\mathbf{Proof~~of~~Proposition}~~\ref{formula for left mutation}~~(ii)~~\mathbf{for~~the~~preinjective~~case}$:

%

\begin{equation*}\begin{split}
&\sum_{j=0}^n(-1)^j(v^{\epsilon(\az)})^{-(n-1)(n-j)}u_{\az}^{+(n-j)}u_{\bz}^-u_{\az}^{+(j)}\\
=&\sum_{j=0}^n(-1)^j(v^{\epsilon(\az)})^{-(n-1)(n-j)+(n-j)(n-j-1)+j(j-1)}u_{(n-j)\az}^+u_{\bz}^-u_{j\az}^+\\
=&\sum_{j=0}^n(-1)^j(v^{\epsilon(\az)})^{j(j-1)-j(n-j)}u_{(n-j)\az}^+u_{\bz}^-u_{j\az}^+\\
=&(\sum_{j=0}^n(-1)^j(v^{\epsilon(\az)})^{j(j-1)}\left|\begin{smallmatrix}
n\\j
\end{smallmatrix}\right])u_{\bz}^-u_{n\az}^++\\
&v^{-\epsilon(\bz)}K_{-\bz}\sum_{j=0}^{n-1}\sum_{i=0}^{n-1-j}(-1)^j(v^{\epsilon(\az)})^{j(j-1)}
\left|\begin{smallmatrix}
i+j\\j
\end{smallmatrix}\right]\sum_{\lz_{n-j,i}:g_{\bz,\lz_{n-j,i}}^{(n-j-i)\az}\neq0}u_{\lambda_{n-j,i}+(i+j)\az}^+.
\end{split}\end{equation*}
By Lemma \ref{zuhe} the first term in the final equation vanishes. Moreover, let $i+j=t$, then using $\sum_{j=0}^{n-1}\sum_{i=0}^{n-1-j}=\sum_{t=0}^{n-1}\sum_{j=0}^t$ we obtain that
\begin{equation*}\begin{split}&\sum_{j=0}^n(-1)^j(v^{\epsilon(\az)})^{-(n-1)(n-j)}u_{\az}^{+(n-j)}u_{\bz}^-u_{\az}^{+(j)}\\
=&v^{-\epsilon(\bz)}K_{-\bz}\sum_{t=0}^{n-1}\sum_{j=0}^{t}(-1)^j(v^{\epsilon(\az)})^{j(j-1)}
\left|\begin{smallmatrix}
t\\j
\end{smallmatrix}\right]\sum_{\lz_{n-j,t-j}:g_{\bz,\lz_{n-j,t-j}}^{(n-t)\az}\neq0}u_{\lambda_{n-j,t-j}+t\az}^+\\
\overset{t=0}{=\!=}& v^{-\epsilon(\bz)}K_{-\bz}u_{\gamma}^+,\end{split}\end{equation*}
where we have used Lemma \ref{zuhe} again for any $t\neq 0$ in the last equation above. This finishes the proof.

\subsubsection{$(V_{\az},V_{\bz})$ lies in the preprojective component}

\begin{Lem}\label{shang}
Let $1\leq j\leq n-1$, and \begin{equation}\label{cok}\xymatrix{0\ar[r]&V_{j\az}\ar[r]^-{\psi}&V_{\bz}\ar[r]^-{\varphi}&V_{c_j}\ar[r]&0}\end{equation} be a short exact sequence in $\modcat \Lambda$. Then
\begin{itemize}
\item[(1)] $\Hom_{\Lz}(V_{\bz},V_{c_j})\cong \End_{\Lz}V_{\bz}\cong\End_{\Lz}V_{c_j};$
\item[(2)] $V_{c_j}$ is indecomposable, $a_{c_j}=a_{\bz}$, and $g_{c_j,j\az}^{\bz}=1.$
\end{itemize}
\end{Lem}

\begin{pf}
Applying the functor $\Hom_{\Lz}(V_{\bz},-)$ to the sequence $(\ref{cok})$, we obtain that
$$\Hom_{\Lz}(V_{\bz},V_{c_j})\cong \End_{\Lz}V_{\bz}.$$
Since $V_{\bz}$ is exceptional, $\End_{\Lz}V_{\bz}$ is a field. So $\Hom_{\Lz}(V_{\bz},V_{c_j})$ is an one-dimensional vector space over $\End_{\Lz}V_{\bz}$. Hence, $V_{c_j}$ is indecomposable, and thus $\End_{\Lz}V_{c_j}$ is a local $k$-algebra. We claim that $\End_{\Lz}V_{c_j}$ is a field. Indeed, for any $0\neq b\in\End_{\Lz}V_{c_j}$, there exists $g\in\End_{\Lz}V_{\bz}$ such that $b\varphi=\varphi g$. Clearly, $g\neq0$, so $g$ is an isomorphism. Suppose that $b$ is nilpotent, say $b^m=0$ for some positive integer $m$. Then $\varphi g^m=b^m \varphi=0$, and thus $\varphi=0$. This is a contradiction.

Applying the functor $\Hom_{\Lz}(-,V_{c_j})$ to the sequence $(\ref{cok})$, we have the following exact sequence:
$$\xymatrix{0\ar[r]&\End_{\Lz}V_{c_j}\ar[r]^-{\varphi^\ast}&\Hom_{\Lz}(V_{\bz},V_{c_j})\ar[r]^-{\psi^\ast}&\Hom_{\Lz}(V_{j\az},V_{c_j})}.$$
We claim that $\psi^\ast=0$. Indeed, for any $f\in\Hom_{\Lz}(V_{\bz},V_{c_j})$, we can write $f=a\varphi$ for some $a\in\End_{\Lz}V_{\bz}$, since $\Hom_{\Lz}(V_{\bz},V_{c_j})$ is one-dimensional over $\End_{\Lz}V_{\bz}$. So $\psi^\ast f=f\psi=a\varphi\psi=0$. Hence, $\End_{\Lz}V_{c_j}\cong\Hom_{\Lz}(V_{\bz},V_{c_j})$. Thus, $\End_{\Lz}V_{c_j}\cong\End_{\Lz}V_{\bz}$ as fields, and then $a_{c_j}=a_{\bz}$.

Since $\Hom_{\Lz}(V_{\bz},V_{c_j})$ as an $\End_{\Lz}V_{c_j}$-vector space is one-dimensional, it is easy to see that $g_{c_j,j\az}^{\bz}=1.$
\end{pf}

\begin{Lem}\label{yin3}
Let $1\leq j\leq n-1$. Then
$$u_{(n-j)\az}^+u_{\bz}^-=u_{\bz}^-u_{(n-j)\az}^++\sum_{i=1}^{n-j}\sum_{c_i:g_{c_i,i\az}^{\bz}\neq0}(v^{\epsilon(\az)})^{-i(n-j)}u_{c_i}^-K_{-i\az}u_{(n-j-i)\az}^+.$$
\end{Lem}

\begin{pf}
Note that $$\Delta(u_{\bz}^-)=u_{\bz}^-\otimes 1+K_{-\bz}\otimes u_{\bz}^-+\sum_{i=1}^{n-j}\sum_{c_i:g_{c_i,i\az}^{\bz}\neq0}b_iu_{c_i}^-K_{-i\az}\otimes u_{i\az}^-+\mbox{else},$$
where $$b_i=v^{\lr{\bz-i\az,i\az}}g_{c_i,i\az}^{\bz}\frac{a_{c_i}a_{i\az}}{a_{\bz}}=(v^{\epsilon(\az)})^{-i^2}a_{i\az}.$$
On the other hand, we have $$\Delta(u_{(n-j)\az}^+)=u_{(n-j)\az}^+\otimes 1+\sum_{i=1}^{n-j-1}a_iu_{(n-j-i)\az}^+K_{i\az}\otimes u_{i\az}^++K_{(n-j)\az}\otimes u_{(n-j)\az}^++\mbox{else},$$
where \begin{equation*}\begin{split}a_i&=(v^{\epsilon(\az)})^{(n-j-i)i}\frac{|\Ext_{\Lz}^1(V_{(n-j-i)\az},V_{i\az})_{V_{(n-j)\az}}|}{
|\Hom_{\Lz}(V_{(n-j-i)\az},V_{i\az})|}\\
&=(v^{\epsilon(\az)})^{(n-j-i)i}/(q^{\epsilon(\az)})^{(n-j-i)i}\\
&=(v^{\epsilon(\az)})^{-(n-j-i)i}.\end{split}\end{equation*}
Hence, by $(\ref{Drinfeld relation})$, we obtain that
\begin{equation*}\begin{split}
u_{(n-j)\az}^+u_{\bz}^-&=u_{\bz}^-u_{(n-j)\az}^++\sum_{i=1}^{n-j}\sum_{c_i:g_{c_i,i\az}^{\bz}\neq0}b_ia_{n-j-i}a_{i\az}^{-1}u_{c_i}^-K_{-i\az}
u_{(n-j-i)\az}^+\\
&=u_{\bz}^-u_{(n-j)\az}^++\sum_{i=1}^{n-j}\sum_{c_i:g_{c_i,i\az}^{\bz}\neq0}(v^{\epsilon(\az)})^{-i^2-(n-j-i)i}u_{c_i}^-K_{-i\az}
u_{(n-j-i)\az}^+\\
&=u_{\bz}^-u_{(n-j)\az}^++\sum_{i=1}^{n-j}\sum_{c_i:g_{c_i,i\az}^{\bz}\neq0}(v^{\epsilon(\az)})^{-i(n-j)}u_{c_i}^-K_{-i\az}
u_{(n-j-i)\az}^+.
\end{split}\end{equation*}
\end{pf}

\begin{Lem}\label{leading term}
$$u_{n\az}^+u_{\bz}^-=u_{\bz}^-u_{n\az}^++\sum_{i=1}^{n-1}\sum_{c_i:g_{c_i,i\az}^{\bz}\neq0}(v^{\epsilon(\az)})^{-ni}u_{c_i}^{-}K_{-i\az}u_{(n-i)\az}^+
+v^{-\epsilon(\bz)}K_{-\bz}u_{\gamma}^+.$$
\end{Lem}

\begin{pf}
Note that $$\Delta(u_{\bz}^-)=u_{\bz}^-\otimes 1+K_{-\bz}\otimes u_{\bz}^-+\sum_{i=1}^{n-1}\sum_{c_i:g_{c_i,i\az}^{\bz}\neq0}b_iu_{c_i}^-K_{-i\az}\otimes u_{i\az}^-+\mbox{else},$$ where $b_i=(v^{\epsilon(\az)})^{-i^2}a_{i\az}$. Moreover,
$$\Delta(u_{n\az}^+)=u_{n\az}^+\otimes 1+\sum_{i=1}^{n-1}a_iu_{(n-i)\az}^+K_{i\az}\otimes u_{i\az}^++K_{n\az}\otimes u_{n\az}^++\tilde{a}_0u_{\bz}^+K_{\gz}\otimes u_{\gz}^++\mbox{else},$$
where $a_i=(v^{\epsilon(\az)})^{-(n-i)i}$, and $\tilde{a}_0=v^{\lr{\bz,n\az-\bz}}\frac{|\Ext_{\Lz}^1(V_{\bz},V_{\gz})_{V_{n\az}}|}{|\Hom_{\Lz}(V_{\bz},V_{\gz})|}
=v^{-\epsilon(\bz)}(q^{\epsilon(\bz)}-1)$.
Hence, by $(\ref{Drinfeld relation})$, we obtain the desired identity.
\end{pf}

\begin{Lem}\label{gexiang}
Let $1\leq j\leq n-1$. Then $u_{(n-j)\az}^+u_{\bz}^-u_{j\az}^+=$
$$(v^{\epsilon(\az)})^{j(n-j)}\left|\begin{smallmatrix}
n\\j
\end{smallmatrix}\right]u_{\bz}^-u_{n\az}^++\sum_{i=1}^{n-j}\sum_{c_i:g_{c_i,i\az}^{\bz}\neq0}
(v^{\epsilon(\az)})^{j(n-j)-ni}\left|\begin{smallmatrix}
n-i\\j
\end{smallmatrix}\right]u_{c_i}^-K_{-i\az}u_{(n-i)\az}^+.$$
\end{Lem}

\begin{pf}
For any positive integers $i$ and $j$, $$u_{i\az}^+u_{j\az}^+=(v^{\epsilon(\az)})^{ij}g_{i\az,j\az}^{(i+j)\az}u_{(i+j)\az}^+
=(v^{\epsilon(\az)})^{ij}\left|\begin{smallmatrix}
i+j\\j
\end{smallmatrix}\right]u_{(i+j)\az}^+.$$
Hence, by Lemma $\ref{yin3}$, we obtain that
\begin{equation*}\begin{split}&u_{(n-j)\az}^+u_{\bz}^-u_{j\az}^+=u_{\bz}^-u_{(n-j)\az}^+u_{j\az}^++\sum_{i=1}^{n-j}\sum_{c_i:g_{c_i,i\az}^{\bz}\neq0}(v^{\epsilon(\az)})^{-i(n-j)}u_{c_i}^-K_{-i\az}u_{(n-j-i)\az}^+u_{j\az}^+\\
=&(v^{\epsilon(\az)})^{j(n-j)}\left|\begin{smallmatrix}
n\\j
\end{smallmatrix}\right]u_{\bz}^-u_{n\az}^++\sum_{i=1}^{n-j}\sum_{c_i:g_{c_i,i\az}^{\bz}\neq0}
(v^{\epsilon(\az)})^{-i(n-j)+(n-j-i)j}\left|\begin{smallmatrix}
n-i\\j
\end{smallmatrix}\right]u_{c_i}^-K_{-i\az}u_{(n-i)\az}^+\\
=&(v^{\epsilon(\az)})^{j(n-j)}\left|\begin{smallmatrix}
n\\j
\end{smallmatrix}\right]u_{\bz}^-u_{n\az}^++\sum_{i=1}^{n-j}\sum_{c_i:g_{c_i,i\az}^{\bz}\neq0}
(v^{\epsilon(\az)})^{j(n-j)-ni}\left|\begin{smallmatrix}
n-i\\j
\end{smallmatrix}\right]u_{c_i}^-K_{-i\az}u_{(n-i)\az}^+.\end{split}\end{equation*}
\end{pf}


$\mathbf{Proof~~of~~Proposition}~~\ref{formula for left mutation}~~(ii)~~\mathbf{for~~the~~preprojective~~case}$:
By Lemmas \ref{leading term} and \ref{gexiang} we get
\begin{equation*}\begin{split}
&\sum_{j=0}^n(-1)^j(v^{\epsilon(\az)})^{-(n-1)(n-j)}u_{\az}^{+(n-j)}u_{\bz}^-u_{\az}^{+(j)}
\\=&\sum_{j=0}^n(-1)^j(v^{\epsilon(\az)})^{j(j-1)-j(n-j)}u_{(n-j)\az}^+u_{\bz}^-u_{j\az}^+\\
=&(\sum_{j=0}^n(-1)^j(v^{\epsilon(\az)})^{j(j-1)}\left|\begin{smallmatrix}
n\\j
\end{smallmatrix}\right])u_{\bz}^-u_{n\az}^++\sum_{i=1}^{n-1}\sum_{c_i:g_{c_i,i\az}^{\bz}\neq0}(v^{\ez(\az)})^{-ni}u_{c_i}^-K_{-i\az}u_{(n-i)\az}^+\\
&+\sum_{j=1}^{n-1}\sum_{i=1}^{n-j}\sum_{c_i:g_{c_i,i\az}^{\bz}\neq0}(-1)^j(v^{\ez(\az)})^{j(j-1)-ni}\left|\begin{smallmatrix}
n-i\\j
\end{smallmatrix}\right]u_{c_i}^-K_{-i\az}u_{(n-i)\az}^++
v^{-\epsilon(\bz)}K_{-\bz}u_{\gz}^+.
\end{split}\end{equation*}
For simplicity we denote the right-hand side of the last equality by $R_1+R_2+R_3+R_4$. By Lemma \ref{zuhe} we have $R_1=0$. In order to finish the proof it remains to show that $R_2+R_3=0$. In fact, by convention, we set $u_{c_n}^-=0$. Then

\begin{equation*}\begin{split}
&R_2+R_3\\
=&\sum_{j=0}^{n-1}\sum_{i=1}^{n-j}\sum_{c_i:g_{c_i,i\az}^{\bz}\neq0}(-1)^j(v^{\ez(\az)})^{j(j-1)-ni}\left|\begin{smallmatrix}
n-i\\j
\end{smallmatrix}\right]u_{c_i}^-K_{-i\az}u_{(n-i)\az}^+  \quad(\text{set\ }t=n-i)\\
=&\sum_{j=0}^{n-1}\sum_{t=j}^{n-1}\sum_{c_{n-t}:g_{c_{n-t},(n-t)\az}^{\bz}\neq0}
(-1)^j(v^{\ez(\az)})^{j(j-1)-n(n-t)}\left|\begin{smallmatrix}
t\\j
\end{smallmatrix}\right]u_{c_{n-t}}^-K_{-(n-t)\az}u_{t\az}^+\\
=&\sum_{t=0}^{n-1}\sum_{j=0}^{t}(-1)^j(v^{\ez(\az)})^{j(j-1)}\left|\begin{smallmatrix}
t\\j\end{smallmatrix}\right]\sum_{c_{n-t}:g_{c_{n-t},(n-t)\az}^{\bz}\neq0}(v^{\ez(\az)})^{-n(n-t)}u_{c_{n-t}}^-K_{-(n-t)\az}u_{t\az}^+\\
=&(v^{\ez(\az)})^{-n^2}u_{c_n}^-K_{-n\az}  \quad(\text{here\ we\ have\ used\ Lemma\ \ref{zuhe}}) \\
=&0.
\end{split}\end{equation*}



\subsection{The case of $0<n(\az,\bz)\Dim V_{\az}<\Dim V_{\bz}$}

That is, $(V_{\az},V_{\bz})$ are the two projective modules in $\modcat \Lambda$ with $V_{\az}$ simple projective. Moreover, $V_{\gz}$ is the simple injective $\Lz$-module which can be determined by an arbitrary short exact sequence of the form
$$0\longrightarrow V_{n\az}\longrightarrow V_{\bz}\longrightarrow V_{\gz}\longrightarrow 0.$$
In this case, Lemma $\ref{shang}$ holds for $1\leq j\leq n$, and Lemmas $\ref{yin3}$ and $\ref{gexiang}$ hold for $0\leq j\leq n-1$. In particular, $V_{c_n}\cong V_{\gz}$. Then the proof of Proposition \ref{formula for left mutation} (iii) is analogous to
the proof of Proposition \ref{formula for left mutation} (ii) for the preprojective case, which will be omitted here.

\begin{Rem} Dually, Proposition \ref{formula for right mutation} can be proved in a quite similar way.
\end{Rem}

\section{The second proof for Proposition \ref{formula for left mutation}}

In this section, we produce a short proof for Proposition \ref{formula for left mutation}, which benefits from the explicit expression of Cramer's isomorphism in Theorem \ref{cramer}.

For the same reason as that in Section 5, we only prove the statements (ii) and (iii) in Proposition \ref{formula for left mutation}.
Thus for the exceptional pair $(V_{\az},V_{\bz})$ in $\modcat\Lz$, we assume $\lr{\az,\bz}>0$ and then $n=n(\az,\bz)$. Denote by $D^b(\sC(V_{\az},V_{\bz}))$ the bounded derived category of $\sC(V_{\az},V_{\bz})$ with the translation functor $[1]$.

\subsection{The case of $n\Dim V_{\az}>\Dim V_{\bz}$}

Set $V_{\bz'}=V_{\bz}[-1]\in D^b(\sC(V_{\az},V_{\bz}))$.
Then $\sC(V_{\az},V_{\bz'})$ is derived equivalent to $\sC(V_{\az},V_{\bz})$.
By Cramer's Theorem \ref{cramer}, we have the isomorphism
$$\Phi:D(\sC(V_{\az},V_{\bz'}))\cong D(\sC(V_{\az},V_{\bz})),$$
which preserves $u_{\az}^{\pm}$ and $u_{\gz}^{\pm}$, and sends $u_{\bz'}^{\pm}$ to $v^{-\epsilon(\bz)}u_{\bz}^{\mp}K_{\pm\bz}.$

Note that $V_{\az}$ and $V_{\bz'}$ are simple in $\sC(V_{\az},V_{\bz'})$, and the left mutation $L(V_{\az},V_{\bz'})=V_{\gz}$, which is given by the standard exact sequence
$$0\longrightarrow V_{\bz'}\longrightarrow V_{\gz}\longrightarrow V_{n\az}\longrightarrow 0.$$
By Proposition \ref{formula for left mutation} (i), the following equality holds in $D(\sC(V_{\az},V_{\bz'}))$:
$$u_{\gz}^{\pm}=\sum_{j=0}^n(-1)^j(v^{\epsilon(\az)})^{(n-j)}u_{\az}^{\pm(n-j)}u_{\bz'}^{\pm}
u_{\az}^{\pm(j)}.$$
Using Cramer's isomorphism $\Phi$, we obtain the following equality in $D(\sC(V_{\az},V_{\bz}))$:
$$u_{\gz}^{\pm}=\sum_{j=0}^n(-1)^j(v^{\epsilon(\az)})^{(n-j)}u_{\az}^{\pm(n-j)}(v^{-\ez(\bz)}u_{\bz}^{\mp}K_{\pm\bz})
u_{\az}^{\pm(j)}.$$
Notice that \begin{equation*}\begin{split}
u_{\az}^{\pm(n-j)}u_{\bz}^{\mp}K_{\pm\bz}
u_{\az}^{\pm(j)}&=v^{(-\bz,(n-j)\az-\bz)}K_{\pm\bz}u_{\az}^{\pm(n-j)}u_{\bz}^{\mp}u_{\az}^{\pm(j)}\\
&=v^{2\ez(\bz)-n(n-j)\ez(\az)}K_{\pm\bz}u_{\az}^{\pm(n-j)}u_{\bz}^{\mp}u_{\az}^{\pm(j)}.
\end{split}\end{equation*}
Hence,  $$u_{\gz}^{\pm}=v^{\epsilon(\bz)}K_{\pm\bz}\sum_{j=0}^n(-1)^j(v^{\epsilon(\az)})^{-(n-1)(n-j)}u_{\az}^{\pm(n-j)}u_{\bz}^{\mp}
u_{\az}^{\pm(j)}.$$

\subsection{The case of $0<n\Dim V_{\az}<\Dim V_{\bz}$}

Set $V_{\az'}=V_{\az}[1]$. Then $\sC(V_{\az'},V_{\bz})$ is derived equivalent to $\sC(V_{\az},V_{\bz})$.
By Cramer's Theorem \ref{cramer}, we have the isomorphism
$$\Phi:D(\sC(V_{\az'},V_{\bz}))\cong D(\sC(V_{\az},V_{\bz})),$$
which preserves $u_{\bz}^{\pm}$ and $u_{\gz}^{\pm}$, and sends $u_{\az'}^{\pm}$ to $v^{\epsilon(\az)}u_{\az}^{\mp}K_{\mp\az}.$

Note that $V_{\az'}$ and $V_{\bz}$ are simple in $\sC(V_{\az'},V_{\bz})$, and the left mutation $L(V_{\az'},V_{\bz})=V_{\gz}$, which is given by the standard exact sequence
$$0\longrightarrow V_{\bz}\longrightarrow V_{\gz}\longrightarrow nV_{\az'}\longrightarrow 0.$$
By Proposition \ref{formula for left mutation} (i), the following equality holds in $D(\sC(V_{\az},V_{\bz'}))$:
$$u_{\gz}^{\pm}=\sum_{j=0}^n(-1)^j(v^{\epsilon(\az)})^{(n-j)}u_{\az'}^{\pm(n-j)}u_{\bz}^{\pm}
u_{\az'}^{\pm(j)},$$ since $\epsilon(\az')=\epsilon(\az).$
Using Cramer's isomorphism $\Phi$, we obtain the following equality in $D(\sC(V_{\az},V_{\bz}))$:
$$u_{\gz}^{\pm}=\sum_{j=0}^n(-1)^j(v^{\epsilon(\az)})^{(n-j)}{(v^{\epsilon(\az)}u_{\az}^{\mp}K_{\mp\az})}^{(n-j)}u_{\bz}^{\pm}
{(v^{\epsilon(\az)}u_{\az}^{\mp}K_{\mp\az})}^{(j)}.$$
Notice that for any positive integer $t$, $$(v^{\ez(\az)}u_{\az}^{\mp}K_{\mp\az})^{(t)}=(v^{\ez(\az)})^{-t^2}K_{\mp t\az}u_{\az}^{\mp(t)}.$$
Thus, $$u_{\gz}^{\pm}=\sum_{j=0}^n(-1)^j(v^{\epsilon(\az)})^{(n-j)-(n-j)^2-j^2} K_{\mp (n-j)\az}u_{\az}^{\mp(n-j)} u_{\bz}^{\pm} K_{\mp j\az}u_{\az}^{\mp(j)}.$$
Moreover, observe that
\begin{equation*}\begin{split}
u_{\az}^{\mp(n-j)} u_{\bz}^{\pm} K_{\mp j\az}=&v^{-(j\az, (n-j)\az-\bz)} K_{\mp j\az} u_{\az}^{\mp(n-j)} u_{\bz}^{\pm}\\
=&(v^{\ez(\az)})^{-2j(n-j)+nj} K_{\mp j\az} u_{\az}^{\mp(n-j)} u_{\bz}^{\pm},
\end{split}\end{equation*}
 and $$(n-j)-(n-j)^2-j^2 -2j(n-j)+nj=-(n-1)(n-j).$$
Hence, $$u_{\gz}^{\pm}=K_{\mp n\az}\sum_{j=0}^n(-1)^j(v^{\epsilon(\az)})^{-(n-1)(n-j)}u_{\az}^{\mp(n-j)}u_{\bz}^{\pm}
u_{\az}^{\mp(j)}.$$

\begin{Rem} Analogously, we can prove Proposition \ref{formula for right mutation} by using Cramer's Theorem \ref{cramer}.
\end{Rem}

\section{The versions of Lie algebras}
Let $\A$ be a finitary hereditary abelian $k$-category as before. Let $D^b(\A)$ be the bounded derived category of $\A$ with the suspension functor $T$. Denote by $\mathcal {R}(\A)$ the root category $D^b(\A)/T^2$. The automorphism $T$ of $D^b(\A)$ induces an automorphism of $\mathcal {R}(\A)$, denoted still by $T$. Then $T^2=1$ in $\mathcal {R}(\A)$.

By using the Ringel--Hall algebra approach, Peng and Xiao \cite{PengXiao2} defined a Lie algebra $\cal{L}(\cR):=\cal{L}(\cR(\A))$ via the root category $\mathcal {R}(\A)$. The nilpotent part of $\cal{L}(\cR)$ has a basis $\{u_{[M]}, u_{[TM]}|M \in\rm{ind}\A\}$.
To simplify the notation, in what follows we write $V_{T\az}:=T(V_{\az})$ and $u_{T\az}:=u_{V_{T\az}}$ for each ${\az}\in\P$.
In this section we will show that the main results in Section 3 also hold under the framework of $\cal{L}(\cR)$.


Let $(V_{\az},V_{\bz})$ be an exceptional pair in $\A$.
Recall that $$n(\az,\bz)=\frac{\lr{\az,\bz}}{\lr{\az,\az}},\ m(\az,\bz)=\frac{\lr{\az,\bz}}{\lr{\bz,\bz}},\ n=|n(\az,\bz)|, \ m=|m(\az,\bz)|.$$
Now, let us recall the left mutation $L_{\az}(\bz)$ and right mutation $R_{\bz}(\az)$ of $(V_{\az},V_{\bz})$ in the root category $\mathcal {R}(\A)$ (see for example \cite{Bondal}).
If $(V_{\az},V_{\bz})$ is orthogonal, i.e., $\Hom(V_{\az},V_{\bz})=0=\Hom(V_{\az},V_{T\bz})$, then $L_{\az}(\bz)=V_{T\bz}$ and $R_{\bz}(\az)=V_{T\az}$. Otherwise, there exists $i=0$ or $1$, such that $\Hom(V_{\az},V_{T^i\bz})\neq 0$. In this case, the left and right mutations are determined respectively by the following canonical triangles
$$\xymatrix{V_{T\bz}\ar[r]& L_{\az}(\bz)\ar[r]& nV_{T^i\az}
\ar[r]& V_{\bz}}$$ and
$$\xymatrix{V_{\az} \ar[r]&  mV_{T^i\bz}
\ar[r]& R_{\bz}(\az)\ar[r]& V_{T\az}.}$$

The following result provides explicit formulas in the Lie algebra $\cal{L}(\cR)$ for the left mutation and right mutation respectively.
For the proof we refer to \cite[Proposition 7.3]{LinPeng} {(see also \cite{R4, PengXiao2})}, where the proof is stated for tubular algebra cases, but it is in fact effective more generally.

\begin{Prop}\label{formula for left mutation--Lie}
Let $\A$ be a finitary hereditary abelian $k$-category. Let $(V_{\az},V_{\bz})$ be an exceptional pair in $\A$ with $\Hom(V_{\az},V_{T^i\bz})\neq 0$ for $i=0$ or 1. Then
$$(-1)^n(n!)u_{[L_{\az}(\bz)]}=({\rm{ad}} u_{T^i\az})^n u_{T\bz}$$ and
$$(m!)u_{[R_{\bz}(\az)]}=({\rm{ad}} u_{T^i\bz})^m u_{T\az}.$$
\end{Prop}

For any exceptional sequence $\cE=(V_{\az_1},V_{\az_2},\cdots,V_{\az_r})$ in $\A$, the subalgebra generated by $\{u_{\az_i}, u_{{T\az_i}}|1\leq i\leq r\}$ is denoted by $\cal{L}_{\cE}(\cR)$.
As an immediate consequence of Proposition \ref{formula for left mutation--Lie}, we have the following

\begin{Prop}\label{Main Prop for Lie algbera}
Let $\A$ be a finitary hereditary abelian $k$-category. If two exceptional sequences $\cE_1$ and $\cE_2$ are mutation equivalent, then $\cal{L}$$_{\cE_1}(\cR)=\cal{L}$$_{\cE_2}(\cR)$.
\end{Prop}

Let $\cal{L}_{\mathscr{E}}(\cR)$ be the subalgebra of $\cal{L}(\cR)$ generated by all $u_{[M]}$'s and
$u_{[TM]}$'s with $M$ exceptional in $\A$, which is called the {\em composition Lie algebra} of $\mathcal {R}(\A)$. By definition, it is easy to see that the composition Lie algebra is invariant under derived equivalences. 
Moreover, in a similar way to the proof in Theorem \ref{add conditions for A}, we obtain the following

\begin{Prop}\label{Lie algebra gen by exc seq}
Assume the conditions in Theorem \ref{add conditions for A} hold. Then for each complete exceptional sequence $\cE$ in $\A$, $L_{\cE}(\cR)=\cal{L}_{\mathscr{E}}(\cR)$.
\end{Prop}

\begin{Rem}
(1) In the case $\A=\modcat A$, the simple $A$-modules form a complete exceptional sequence when suitably ordered. Hence the composition Lie algebra can be defined via simple modules, {which in fact has been adopted as its original definition, see for example \cite{Ringel1990, PengXiao2}.}

(2) In the case $\A=\coh\bbX$, the composition Lie algebra is a nice model to realize the loop algebra of the Kac--Moody algebra associated to $\bbX$, for details we refer to Subsection \ref{composition lie algebra subsection}.
\end{Rem}

\section{Applications}

\subsection{Analogues of quantum Serre relations }

We provide some analogues of quantum Serre relations for any exceptional pair $(V_{\az},V_{\bz})$ in $\A$.
\begin{Prop}\label{Serre's relation}
Let $(V_{\az},V_{\bz})$ be an exceptional pair in $\A$,
\begin{itemize}
\item[(i)] if $\lr{\az,\bz}\leq0$, then
$$\sum_{j=0}^{n+1}(-1)^ju_{\az}^{\pm(n+1-j)}u_{\bz}^{\pm}u_{\az}^{\pm(j)}=0
=\sum_{j=0}^{m+1}(-1)^ju_{\bz}^{\pm(m+1-j)}u_{\az}^{\pm}u_{\bz}^{\pm(j)};$$
\item[(ii)] if $\lr{\az,\bz}>0$, then
$$\sum_{j=0}^{n+1}(-1)^j(v^{\ez({\az})})^{nj}u_{\az}^{\pm(n+1-j)}u_{\bz}^{\mp}u_{\az}^{\pm(j)}=0
=\sum_{j=0}^{m+1}(-1)^j(v^{\ez({\bz})})^{-mj}u_{\bz}^{\pm(m+1-j)}u_{\az}^{\mp}u_{\bz}^{\pm(j)}.$$
\end{itemize}
\end{Prop}

\begin{pf}
If $\lr{\az,\bz}\leq0$, then $(V_{\az},V_{\bz})$ are two relative simple objects in $\sC(V_{\az},V_{\bz})$, thus the first statement is well-known (see \cite{Ringel1,Ringel2}). If $\lr{\az,\bz}>0$, then one can use a similar proof as that for Proposition \ref{formula for left mutation} in Section 5 or Section 6 to obtain the result.
\end{pf}

\subsection{Lusztig's symmetries}
First of all, let us recall the isomorphism between the double composition algebra $D_{\mathscr{E}}(A):=D_{\mathscr{E}}(\modcat A)$ and the corresponding quantum group $U_q(\mathfrak{g})$ (cf. \cite{X97}), which is defined on generators by
$$\theta:D_{\mathscr{E}}(A)\longrightarrow U_q({\fg}),~~u_i^+\mapsto E_i,~~u_i^-\mapsto -v^{-\epsilon(i)}F_i,~~K_i\mapsto \tilde{K}_i.$$
The notations used here for elements of $U_q(\mathfrak{g})$ are the same as those in \cite[Chapter 3]{Lus}. Lusztig \cite[Chapter 37]{Lus} defined four families of symmetries as automorphisms of $U_q(\mathfrak{g})$, which are denoted by $T_{i,e}^{'}$ and $T_{i,e}^{''}$, where $e=\pm1$ and $i\in I$.

In the left mutation formula in Proposition $\ref{formula for left mutation}$(i), $V_{\az}$ and $V_{\bz}$ are two simple $\Lambda$-modules, say $\az=i$ and $\bz=j$. By the standard exact sequence
$$0\longrightarrow V_j\longrightarrow V_{\gz}\longrightarrow nV_i\longrightarrow 0,$$ where
$n=-\frac{2(i,j)}{(i,i)}$, we obtain that
$\dim_kV_{\gz}=\dim_kV_j+n\dim_kV_i=\epsilon(j)+n\epsilon(i)$. Besides,
$\epsilon(\gz)=\lr{\gz,\gz}=\lr{j+ni,j+ni}=\epsilon(j)+n\cdot(i,j)+n^2\epsilon(i)=\epsilon(j)
$. So $-\dim_kV_{\gz}+\epsilon(\gz)=-n\epsilon(i)$. Hence,
\begin{equation*}\begin{split}\theta(u_{\gz}^-)&=\sum_{r=0}^n(-1)^r(v^{\epsilon(i)})^{n-r}(-v^{-\epsilon(i)}F_i)^{(n-r)}(-v^{-\epsilon(j)}F_j)
(-v^{-\epsilon(i)}F_i)^{(r)}\\&=
(-1)^{n+1}v^{-\epsilon(j)}\sum_{r=0}^n(-1)^rv^{-r\epsilon(i)}F_i^{(n-r)}F_jF_i^{(r)}\\&=
(-1)^{n+1}v^{-\epsilon(j)}T_{i,1}^{'}(F_j).\end{split}\end{equation*}
Since $\epsilon(j)=\epsilon(\gz)$, we obtain that
\begin{equation}\label{L1}T_{i,1}^{'}(F_j)=(-1)^{n+1}v^{\epsilon(\gz)}\theta(u_{\gz}^-).\end{equation}
Similarly, we obtain that (see also \cite{CX,R4})
\begin{equation}T_{i,1}^{''}(E_j)=v^{-\dim_kV_{\gz}+\epsilon(\gz)}\theta(u_{\gz}^+).\end{equation}

In the right mutation formula in Proposition $\ref{formula for right mutation}$(i), $V_{\az}$ and $V_{\bz}$ are two simple $\Lambda$-modules, say $\az=i$ and $\bz=j$. By the standard exact sequence
$$0\longrightarrow mV_j\longrightarrow V_{\lz}\longrightarrow V_i\longrightarrow 0,$$ where
$m=-\frac{2(i,j)}{(j,j)}$, we obtain that $\dim_kV_{\lz}=\epsilon(i)+m\epsilon(j)$ and $\epsilon(\lz)=\epsilon(i)$, thus $-\dim_kV_{\lz}+\epsilon(\lz)=-m\epsilon(j)$. Hence,
\begin{equation*}\begin{split}\theta(u_{\lz}^-)&=\sum_{r=0}^m(-1)^r(v^{\epsilon(j)})^{m-r}(-v^{-\epsilon(j)}F_j)^{(r)}(-v^{-\epsilon(i)}F_i)
(-v^{-\epsilon(j)}F_j)^{(m-r)}\\
&=(-1)^{m+1}v^{-\epsilon(i)}\sum_{r=0}^m(-1)^rv^{-r\epsilon(j)}F_j^{(r)}F_iF_j^{(m-r)}\\&=
(-1)^{m+1}v^{-\epsilon(i)}T_{i,-1}^{''}(F_j).\end{split}\end{equation*}
Since $\epsilon(i)=\epsilon(\lz)$, we obtain that \begin{equation}T_{i,-1}^{''}(F_j)=(-1)^{m+1}v^{\epsilon(\lz)}\theta(u_{\lz}^-).\end{equation}
Similarly, we obtain that (see also \cite{CX,R4})
\begin{equation}\label{L4}T_{i,-1}^{'}(E_j)=v^{-\dim_kV_{\lz}+\epsilon(\lz)}\theta(u_{\lz}^+).\end{equation}
We remark that all four kinds of Lusztig's symmetries have appeared in Equations $(\ref{L1}-\ref{L4})$.

\subsection{Double composition algebras for weighted projective lines}\label{composition algebra of coh subsection}


The \emph{composition algebra} $C(\bbX)$ of $\coh\bbX$ is defined by Schiffmann \cite{Sch2004} as the subalgebra of the Ringel--Hall algebra $\cH(\coh\bbX)$ generated by $u_{\co(l\vc)}, l\in \bbZ$; $T_r, r\in\bbZ_{>0}$ and $u_{S_{ij}}, 1\leq i\leq t, 1\leq j\leq p_i-1$, which has been used to realize a certain ``positive part" of the quantum loop algebra of the Kac-Moody algebra associated to $\bbX$.

Later on, Burban and Schiffmann \cite{BurSch} defined an extended version of the composition algebra, denoted by {$U(\bbX):=C(\coh\bbX)\otimes \bf{K}$}, where $\bfK$ is the group algebra $\bbC[K(\coh\bbX)]$. They showed that $U(\bbX)$ is a topological bialgebra and then defined $DU(\bbX)$ as its \emph{reduced Drinfeld double}. Meanwhile, Dou, Jiang and Xiao \cite{DJX} defined the \emph{``double composition algebra"} $DC(\bbX)$ of $\coh\bbX$ as the subalgebra of $D(\coh\bbX)$ generated by two copies of $C(\bbX)$, say $C^{\pm}(\bbX)$, together with the group algebra $\bfK$.


\begin{Prop}\label{identity for composition algebra for wpl}
Keep notations as above. The following subalgebras of $D(\coh\bbX)$ coincide: $DC(\bbX)=DU(\bbX)=D_{\mathscr{E}}(\coh\bbX)$.
\end{Prop}

\begin{pf} It has been mentioned in \cite{DJX} that $DC(\bbX)=DU(\bbX)$. By \cite[Corollary 5.23]{BurSch}, the subalgebra $DU(\bbX)$ is generated by $u^{\pm}_{\co(\vx)}, \vx\in\bbL$ together with the torus {$\bf{K}$}. Note that each line bundle $\co(\vx)$ is exceptional in $\coh\bbX$. Hence by the definition of $D_{\mathscr{E}}(\coh\bbX)$, it contains $DU(\bbX)$ as a subalgebra. On the other hand, since $\bigoplus_{0\leq\vx\leq \vc}\co(\vx)$ forms a tilting object in $\coh\bbX$, 
they can be arranged as a complete exceptional sequence of $\coh\bbX$. Then by the Main Theorem we obtain that $D_{\mathscr{E}}(\coh\bbX)\subseteq DU(\bbX)$, which completes the proof.
\end{pf}

\begin{Rem}
{The above proposition shows that for the category $\coh\bbX$, the reduced Drinfeld double composition algebra $DU(\bbX)$ defined by Burban and Schiffmann, and the ``double composition algebra" $DC(\bbX)$ defined by Dou, Jiang and Xiao, both coincide with the double composition algebra $D_{\mathscr{E}}(\coh\bbX)$ defined via exceptional objects. Moreover, by Corollary \ref{main examples}, they can be generated by any complete exceptional sequence in $\coh\bbX$.}
\end{Rem}

\subsection{Composition Lie algebras for weighted projective lines}\label{composition lie algebra subsection}
Now we consider the composition Lie algebra $\cal{L}(\cR):=\cal{L}(\cR(\coh\bbX))$ of $\coh\bbX$. In Theorem 2 of \cite{CB}, Crawley-Boevey found a series of elements in $\cal{L}(\cR)$ to satisfy the generating relations of the loop algebra of the Kac--Moody algebra associated to $\bbX$.
We emphasize that the subalgebra of $\cal{L}(\cR)$ generated by these elements, which we denote by $\cal{L}$$_c(\cR)$, is nothing else but the composition Lie algebra $\cal{L}_{\mathscr{E}}(\cR)$.

\begin{Prop}\label{composition Lie algebra for wpl}
Keep notations as above. Then $\cal{L}$$_c(\cR)=\cal{L}_{\mathscr{E}}(\cR)$.
\end{Prop}

\begin{pf}
We first claim that $\cal{L}_{\mathscr{E}}(\cR)\subseteq\cal{L}$$_c(\cR)$. In fact, note that $\cal{L}$$_c(\cR)$ contains the elements $u_{[M]}$ and $u_{[TM]}$ for $M=\co, \co(\vc)$ or $S_{ij}, 1\leq i\leq t, 1\leq j\leq p_i-1$. Moreover, these elements form a complete exceptional sequence in the following way
$$(\co,\co(\vc);S_{1,p_1-1},\cdots, S_{11};\cdots; S_{t,p_t-1},\cdots, S_{t1}).$$ Hence the claim follows from Proposition \ref{Lie algebra gen by exc seq}.

On the other hand, we need to show that $\cal{L}$$_c(\cR)\subseteq \cal{L}_{\mathscr{E}}(\cR).$ By definition we know that $c, h_{\star 0}, e_{\star r}, f_{\star r}$ and  $e_{v 0}, f_{v0}$ belong to the $\cal{L}_{\mathscr{E}}(\cR)$ for $r\in\bbZ, v=ij, 1\leq i\leq t, 1\leq j\leq p_i-1$. It suffices to show that the other generators of $\cal{L}$$_c(\cR)$ can be generated by these elements. In fact, for any $r\neq 0$, $h_{\star r}=[e_{\star  r}, f_{\star  0}]\in\cal{L}_{\mathscr{E}}(\cR)$. Moreover, by induction on the index $j$ (we write $\star=i0$ for convenience), we can obtain $e_{v r}, f_{v r}, h_{v r}$ belong to $\cal{L}_{\mathscr{E}}(\cR)$ for $v=ij$ and $r\neq 0$ by using the relations $e_{vr}=[e_{v0}, h_{v'r}], f_{vr}=[h_{v'r}, f_{v0}]\ (v'=i,j-1)$ and $h_{vr}=[e_{vr}, f_{v0}]$.
\end{pf}

\subsection{Bridgeland's Hall algebras and Modified Ringel--Hall algebras}
Let $\A$ be a finitary hereditary abelian $k$-category with enough projectives. In order to give an intrinsic realization of the entire quantized enveloping algebra via Hall algebras, Bridgeland \cite{Br} considered the Hall algebra of 2-cyclic complexes of projective objects in $\A$, and achieved an algebra $\mathcal {D}\mathcal {H}_{{\rm red}}(\mathcal{A})$, called the {\em reduced Bridgeland's Hall algebra} of $\A$, by taking some localization and reduction. By \cite{Yan,ZHC}, $\mathcal {D}\mathcal {H}_{{\rm red}}(\mathcal{A})$ is isomorphic to the reduced Drinfeld double Hall algebra $D(\A)$. In order to generalize Bridgeland's construction to any hereditary abelian category $\A$, Lu and Peng \cite{LP} defined an algebra $\cM\cH_{\bbZ/2,{\rm tw,red}}(\A)$, called the {\em reduced modified Ringel--Hall algebra} of $\A$, and proved that it is also isomorphic to $D(\A)$.
Hence, our results can also be applied to $\mathcal {D}\mathcal {H}_{{\rm red}}(\mathcal{A})$ and $\cM\cH_{\bbZ/2,{\rm tw,red}}(\A)$.

\section*{Acknowledgments}

The authors are grateful to Bangming Deng, Jie Xiao and Jie Sheng for their helpful comments. In particular, Section 6 comes from Deng's suggestion.


\end{document}